\documentclass[12pt]{article}
\usepackage{amsfonts}
\usepackage{amssymb}
\usepackage{myart}

\normalbaselines
\input tcilatex
\begin{document}

\title{G\"{o}del's theorem as a corollary of impossibility of complete
axiomatization of geometry}
\author{Yuri A.Rylov}
\date{Institute for Problems in Mechanics, Russian Academy of Sciences,\\
101-1, Vernadskii Ave., Moscow, 119526, Russia.\\
e-mail: rylov@ipmnet.ru\\
Web site: {$http://rsfq1.physics.sunysb.edu/\symbol{126}rylov/yrylov.htm$}\\
or mirror Web site: {$http://gasdyn-ipm.ipmnet.ru/\symbol{126}%
rylov/yrylov.htm$}}
\maketitle

\begin{abstract}
Not any geometry can be axiomatized. The paradoxical Godel's theorem starts
from the supposition that any geometry can be axiomatized and goes to the
result, that not any geometry can be axiomatized. One considers example of
two close geometries (Riemannian geometry and $\sigma $-Riemannian one),
which are constructed by different methods and distinguish in some details.
The Riemannian geometry reminds such a geometry, which is only a part of the
full geometry. Such a possibility is covered by the Godel's theorem.
\end{abstract}

\section{Introduction}

Let there be a set $\Omega $ of points $P_{1},P_{2},...$ A geometrical
object $O$ is a subset $O$ of points $P_{1},P_{2},...$ of $\Omega $, $%
(O\subset \Omega )$.

\begin{definition}
Geometry $\mathcal{G}$ is an infinite set $\mathcal{S}_{\mathcal{G}}$ of
prepositions $\mathcal{P}_{1},\mathcal{P}_{2},...$ on properties of
geometrical objects $O_{1},O_{2},...\subset \Omega $.
\end{definition}

\begin{definition}
If one can choose a finite or countable set $\mathcal{S}_{\mathcal{A}}%
\mathcal{\subset S}_{\mathcal{G}}$ of prepositions in such a way, that the
infinite set $\mathcal{S}_{\mathcal{G}}$ of all prepositions $\mathcal{P}%
_{1},\mathcal{P}_{2},...$ on properties of geometrical objects $%
O_{1},O_{2},...\subset \Omega $ can be obtained as a result of logical
reasonings of prepositions $\mathcal{P}\in \mathcal{S}_{\mathcal{A}}$, the
geometry $\mathcal{G}$ may be axiomatized, and the set $\mathcal{S}_{%
\mathcal{A}}$ of prepositions $\mathcal{P}$ forms axiomatics (the system of
axioms) of the geometry $\mathcal{G}$.
\end{definition}

It is evident, that not any geometry can be axiomatized. However, there are
geometries, which can be axiomatized. For instance, the proper Euclidean
geometry can be axiomatized.

Let us consider a geometry $\mathcal{G}$, which can be axiomatized partly.
It means that a part $\mathcal{S}_{\mathcal{G}_{1}}\subset \mathcal{S}_{%
\mathcal{G}}$ of all propositions $\mathcal{P}$ of geometry $\mathcal{G}$
may be obtained as a result of logical reasonings of propositions $\mathcal{P%
}\in \mathcal{S}_{\mathcal{A}}$, where $\mathcal{S}_{\mathcal{A}}$ is the
set of axioms $\mathcal{S}_{\mathcal{A}}\mathcal{\subset S}_{\mathcal{G}}$.
The set of axioms $\mathcal{S}_{\mathcal{A}}$ is supposed to be complete in
the sense, that any supposition $\mathcal{P}\in \mathcal{S}_{\mathcal{G}_{1}}
$ may be deduced from the set $\mathcal{S}_{\mathcal{A}}$ of axioms. It is
supposed that the set $\mathcal{S}_{\mathrm{r}}=\mathcal{S}_{\mathcal{G}%
}\backslash \mathcal{S}_{\mathcal{G}_{1}}$ of remaining propositions of
geometry $\mathcal{G}$ cannot be deduced from the set $\mathcal{S}_{\mathcal{%
A}}$ of axioms.

The Godel's theorem proves the same, but Godel starts from the supposition,
that the complete axiomatization of geometry is possible, and he deduced
essentially that the complete axiomatization is not always possible.

The real problem of the geometry construction lies in the fact, that we know
no other method of the geometry construction except for the deduction of the
geometry propositions from the geometry axioms. Euclid had deduced his
Euclidean geometry from the system of axioms. It was proved \cite{H30}, that
the Euclidean axioms are compatible between themselves, and the proper
Euclidean geometry is a consistent geometry.

Mathematicians try to repeat the Euclidean experience for construction of
other (non-uniform) geometries. They use the Euclidean method of the
geometry construction. As a result one cannot be sure that one deduces the
full geometry, but not only its part.

We start from the general supposition that the complete axiomatization of a
geometry is not always possible. In this case the result of the Godel's
theorem is evident. However, to adduce such a supposition, one is to have an
alternative method of the geometry construction, because a formal
consideration of the set of the geometrical propositions is not constructive
without a method of these propositions obtaining.

In reality there is an alternative method of the geometry construction. It
is based on two statements.

\begin{enumerate}
\item Geometry is described completely by its metric (distance between any
two points of the space).

\item The proper Euclidean geometry is a true geometry, which can be
described completely by its metric.
\end{enumerate}

We shall refer to geometry, which is completely described by its metric, as
the physical geometry. The geometry, which is deduced from the system of
some axioms, will be referred to as an axiomatic geometry, or a mathematical
geometry. Then the alternative method of the physical geometry construction
is formulated as follows.

\textit{Any physical geometry is obtained from the proper Euclidean geometry
as a result of its deformation.}

It means that all propositions of the proper Euclidean geometry are
expressed via Euclidean metric, and the Euclidean metric in all propositions
is replaced by the metric of the geometry in question. There is a theorem,
which proves that all propositions of the proper Euclidean geometry can be
expressed via Euclidean metric \cite{R2001}. As a result one obtains all
propositions of the geometry in question, expressed via its metric.

Usually one supposes, that the term "metric" means the distance, satisfying
the triangle axiom. In the previous presentation the term "metric" means the
distance simply (which is not satisfies the triangle axiom, in general).

To avoid a confusion, one uses the world function $\sigma \left( P,Q\right) =%
\frac{1}{2}\rho ^{2}\left( P,Q\right) $ instead of the metric $\rho \left(
P,Q\right) $. It is more convenient from technical point of view. Besides,
in the case of indefinite geometry (for instance, in geometry of Minkowski)
the world function $\sigma $ is real, although $\rho =\sqrt{2\sigma }$ may
be imaginary.

\section{Deformation principle as a method of a physical geometry
construction}

If one knows the proper Euclidean geometry, the expression of the proper
Euclidean propositions via the Euclidean world function is a purely
technical problem. But there are some subtleties in this problem. The fact
is that, the proper Euclidean geometry has specific Euclidean properties and
general geometric properties. All propositions of the proper Euclidean
geometry are to be expressed only via general geometric properties. Only in
this case the expression via Euclidean world function may be deformed and
used in the deformed geometry. The fact is that, the specific properties of
the proper Euclidean geometry are different in the Euclidean spaces of
different dimensions. Any specific property of the proper Euclidean geometry
contains a reference to the dimension $n$ of the Euclidean space.

The general geometric propositions of the proper Euclidean geometry do not
refer to the dimension $n$ of the Euclidean space. Thus, only the Euclidean
propositions, which do not contain a reference to the dimension $n$, may be
deformed to obtain corresponding relation of the deformed geometry.

We consider a simple example. Vector $\mathbf{PQ}=\overrightarrow{PQ}$ is an
ordered set of two points $P$ and $Q$. The scalar product $\left( \mathbf{P}%
_{0}\mathbf{P}_{1}.\mathbf{Q}_{0}\mathbf{Q}_{1}\right) $ of two vectors $%
\mathbf{P}_{0}\mathbf{P}_{1}$ and $\mathbf{Q}_{0}\mathbf{Q}_{1}$ is defined
by the relation%
\begin{equation}
\left( \mathbf{P}_{0}\mathbf{P}_{1}.\mathbf{Q}_{0}\mathbf{Q}_{1}\right)
=\sigma \left( P_{0},Q_{1}\right) +\sigma \left( P_{1},Q_{0}\right) -\sigma
\left( P_{0},Q_{0}\right) -\sigma \left( P_{1},Q_{1}\right)  \label{a1.1}
\end{equation}%
where $\sigma $ is the world function%
\begin{equation}
\sigma :\qquad \Omega \times \Omega \rightarrow \mathbb{R},\qquad \sigma
\left( P,Q\right) =\sigma \left( Q,P\right) ,\qquad \sigma \left( P,P\right)
=0,\qquad \forall P,Q\in \Omega  \label{a1.2}
\end{equation}%
Definition (\ref{a1.1}) of the scalar product of two vectors coincides with
the conventional scalar product of vectors in the proper Euclidean space.
(One can verify this easily). The relation (\ref{a1.1}) does not contain a
reference to the dimension of the Euclidean space and to a coordinate system
in it. Hence, the relation (\ref{a1.1}) is a general geometric relation,
which may be considered as a definition of the scalar product of two vectors
in any physical geometry.

Equivalence (equality) of two vectors $\mathbf{P}_{0}\mathbf{P}_{1}$ and $%
\mathbf{Q}_{0}\mathbf{Q}_{1}$ is defined by the relations%
\begin{equation}
\mathbf{P}_{0}\mathbf{P}_{1}\text{eqv}\mathbf{Q}_{0}\mathbf{Q}_{1}:\qquad
\left( \mathbf{P}_{0}\mathbf{P}_{1}.\mathbf{Q}_{0}\mathbf{Q}_{1}\right)
=\left\vert \mathbf{P}_{0}\mathbf{P}_{1}\right\vert \cdot \left\vert \mathbf{%
Q}_{0}\mathbf{Q}_{1}\right\vert \wedge \left\vert \mathbf{P}_{0}\mathbf{P}%
_{1}\right\vert =\left\vert \mathbf{Q}_{0}\mathbf{Q}_{1}\right\vert
\label{a1.3}
\end{equation}%
where $\left\vert \mathbf{P}_{0}\mathbf{P}_{1}\right\vert $ is the length of
the vector $\mathbf{P}_{0}\mathbf{P}_{1}$%
\begin{equation}
\left\vert \mathbf{P}_{0}\mathbf{P}_{1}\right\vert =\sqrt{\left( \mathbf{P}%
_{0}\mathbf{P}_{1}.\mathbf{P}_{0}\mathbf{P}_{1}\right) }=\sqrt{2\sigma
\left( P_{0},P_{1}\right) }  \label{a1.4}
\end{equation}

In the developed form the condition (\ref{a1.3}) of equivalence of two
vectors $\mathbf{P}_{0}\mathbf{P}_{1}$ and $\mathbf{Q}_{0}\mathbf{Q}_{1}$
has the form%
\begin{eqnarray}
\sigma \left( P_{0},Q_{1}\right) +\sigma \left( P_{1},Q_{0}\right) -\sigma
\left( P_{0},Q_{0}\right) -\sigma \left( P_{1},Q_{1}\right) &=&2\sigma
\left( P_{0},P_{1}\right)  \label{a1.5} \\
\sigma \left( P_{0},P_{1}\right) &=&\sigma \left( Q_{0},Q_{1}\right)
\label{a1.6}
\end{eqnarray}

If the points $P_{0},P_{1}$, determining the vector $\mathbf{P}_{0}\mathbf{P}%
_{1}$, and the origin $Q_{0}$ of the vector $\mathbf{Q}_{0}\mathbf{Q}_{1}$
are given, we can determine the vector $\mathbf{Q}_{0}\mathbf{Q}_{1}$, which
is equivalent (equal) to the vector $\mathbf{P}_{0}\mathbf{P}_{1}$, solving
two equations (\ref{a1.5}), (\ref{a1.6}) with respect to the position of the
point $Q_{1}$.

In the case of the proper Euclidean space there is one and only one solution
of equations (\ref{a1.5}), (\ref{a1.6}) independently of its dimension $n$.
In the case of arbitrary physical geometry one can guarantee neither
existence nor uniqueness of the solution of equations (\ref{a1.5}), (\ref%
{a1.6}) for the point $Q_{1}$. Number of solutions depends on the form of
the world function $\sigma $. This fact means a multivariance of the
property of two vectors equivalence in the arbitrary physical geometry. In
other words, single-variance of the vector equality in the proper Euclidean
space is a specific property of the proper Euclidean geometry, and this
property is conditioned by the form of the Euclidean world function. In
other physical geometries this property does not take place, in general.

The multivariance is a general property of a physical geometry. It is
connected with a necessity of solution of algebraic equations, containing
the world function. As far as the world function is different in different
physical geometries, the solution of these equations may be not unique, or
it may not exist at all.

In the proper Euclidean geometry the equality of two vectors $\mathbf{P}_{0}%
\mathbf{P}_{1}$ and $\mathbf{Q}_{0}\mathbf{Q}_{1}$ may be defined also as
equality of components of these vectors at some rectilinear coordinate
system. (It is the conventional method of definition of two vectors
equality).

Let the dimension of the Euclidean space be equal to $n$. Let us introduce $%
n $ linear independent vectors $\mathbf{P}_{0}\mathbf{P}_{k}$, $k=1,2,...n$.
Linear independence of vectors $\mathbf{P}_{0}\mathbf{P}_{k}$ means that the
Gram's determinant 
\begin{equation}
\det \left\vert |\left( \mathbf{P}_{0}\mathbf{P}_{i}.\mathbf{P}_{0}\mathbf{P}%
_{k}\right) |\right\vert \neq 0,\qquad i,k=1,2,...n  \label{a1.7}
\end{equation}%
We construct rectilinear coordinate system with basic vectors $\mathbf{P}_{0}%
\mathbf{P}_{k}$, $k=1,2,...n$ in the $n$-dimensional Euclidean space.
Covariant coordinates $x_{k}=\left( \mathbf{P}_{0}\mathbf{P}_{1}\right) _{k}$
and $y_{k}=\left( \mathbf{Q}_{0}\mathbf{Q}_{1}\right) _{k}$ of vectors $%
\mathbf{P}_{0}\mathbf{P}_{1}$ and $\mathbf{Q}_{0}\mathbf{Q}_{1}$ in this
coordinate system have the form%
\begin{equation}
x_{k}=\left( \mathbf{P}_{0}\mathbf{P}_{1}\right) _{k}=\left( \mathbf{P}_{0}%
\mathbf{P}_{1}.\mathbf{P}_{0}\mathbf{P}_{k}\right) ,\qquad y_{k}=\left( 
\mathbf{Q}_{0}\mathbf{Q}_{1}\right) _{k}=\left( \mathbf{Q}_{0}\mathbf{Q}_{1}.%
\mathbf{P}_{0}\mathbf{P}_{k}\right) ,\qquad k=1,2,...n  \label{a1.8}
\end{equation}%
Equality of vectors $\mathbf{P}_{0}\mathbf{P}_{1}$ and $\mathbf{Q}_{0}%
\mathbf{Q}_{1}$ is written in the form%
\begin{equation}
\left( \mathbf{P}_{0}\mathbf{P}_{1}.\mathbf{P}_{0}\mathbf{P}_{k}\right)
=\left( \mathbf{Q}_{0}\mathbf{Q}_{1}.\mathbf{P}_{0}\mathbf{P}_{k}\right)
,\qquad k=1,2,...n  \label{a1.9}
\end{equation}%
and according to (\ref{a1.1}) it may be written in terms of the world
function (metric). The points $P_{0},P_{1},Q_{0}$ are supposed to be given.
The point $Q_{1}$ is to be determined by equations (\ref{a1.9}). Equations (%
\ref{a1.9}) determine equality of vectors $\mathbf{P}_{0}\mathbf{P}_{1}$ and 
$\mathbf{Q}_{0}\mathbf{Q}_{1}$ only in the $n$-dimensional Euclidean space $%
E_{n}$. Already in the $(n+1)$-dimensional Euclidean space $E_{n+1}$ $n$
equations (\ref{a1.9}) do not determine equality of vectors $\mathbf{P}_{0}%
\mathbf{P}_{1}$ and $\mathbf{Q}_{0}\mathbf{Q}_{1}$, because in $(n+1)$%
-dimensional Euclidean space $E_{n+1}$ one needs $(n+1)$ equations of the
form (\ref{a1.9}) to define equality of vectors. There may be no dimension
in the physical geometry, or the dimension of the space may be different at
different points. In this cases conventional conditions (\ref{a1.9}) of two
vectors equality cannot be used also.

From formal point of view the equations (\ref{a1.9}) define some geometrical
object, or a set of geometrical objects, whose points are described by means
of running point $Q_{1}$. This geometrical object depends on parameters $%
Q_{0},P_{0},P_{1},...P_{n}$. How can one interpret this object? It is quite
unclear.

From formal viewpoint the relations (\ref{a1.5}), (\ref{a1.6}) describe some
geometrical object by means of the running point $Q_{1}$. This object
depends on parameters $Q_{0},P_{0},P_{1}$, and one interprets this as a set
of vectors $\mathbf{Q}_{0}\mathbf{Q}_{1}$, which are equivalent to vector $%
\mathbf{P}_{0}\mathbf{P}_{1}$. One may consider, that the relations (\ref%
{a1.5}), (\ref{a1.6}) describe some geometrical object by means of the
running point $P_{1}$. This object depends on parameters $Q_{0},Q_{1},P_{0}$%
, and one interprets this as a set of vectors $\mathbf{P}_{0}\mathbf{P}_{1}$%
, which are equivalent to vector $\mathbf{Q}_{0}\mathbf{Q}_{1}$.

It is to note that the proper Euclidean geometry is a degenerate geometry in
the sense that the same geometrical object may be described by different
ways in terms of the world function. For instance, the circular cylinder $%
\mathcal{CY}\left( P_{0},P_{1},Q\right) $ is defined by the relation%
\begin{equation}
\mathcal{CY}\left( P_{0},P_{1},Q\right) =\left\{
R|S_{P_{0}P_{1}R}=S_{P_{0}P_{1}Q}\right\}  \label{a1.10}
\end{equation}%
where $P_{0},P_{1}$ are two different points on the axis of the cylinder and 
$Q$ is some point on the surface of the cylinder. The quantity $%
S_{P_{0}P_{1}Q}$ is the area of the triangle with the vertices at the points 
$P_{0},P_{1},Q$. The area of triangle may be calculated by means of the
Hero's, expressing the triangle area via length of the triangle sides, or by
the formula%
\begin{equation}
S_{P_{0}P_{1}Q}=\frac{1}{2}\sqrt{\left\vert 
\begin{array}{ll}
\left( \mathbf{P}_{0}\mathbf{P}_{1}.\mathbf{P}_{0}\mathbf{P}_{1}\right) & 
\left( \mathbf{P}_{0}\mathbf{P}_{1}.\mathbf{P}_{0}\mathbf{Q}\right) \\ 
\left( \mathbf{P}_{0}\mathbf{Q}.\mathbf{P}_{0}\mathbf{P}_{1}\right) & \left( 
\mathbf{P}_{0}\mathbf{Q}.\mathbf{P}_{0}\mathbf{Q}\right)%
\end{array}%
\right\vert }  \label{a1.11}
\end{equation}%
which may be expressed in terms of the world function by means of (\ref{a1.1}%
). In the proper Euclidean geometry the circular cylinder $\mathcal{CY}%
\left( P_{0},P_{1},Q\right) $ depends only on its axis $\mathcal{T}%
_{P_{0}P_{1}}$, passing through the points $P_{0},P_{1}$, but not on
positions of the points $P_{0},P_{1}$ on the axis $\mathcal{T}_{P_{0}P_{1}}$.

The axis of the cylinder $\mathcal{T}_{P_{0}P_{1}}$ is described by the
relation 
\begin{equation}
\mathcal{T}_{P_{0}P_{1}}=\left\{ R|\mathbf{P}_{0}\mathbf{P}_{1}\parallel 
\mathbf{P}_{0}\mathbf{R}\right\}  \label{a1.12}
\end{equation}%
where $\mathbf{P}_{0}\mathbf{P}_{1}\parallel \mathbf{P}_{0}\mathbf{R}$ means
that the vectors $\mathbf{P}_{0}\mathbf{P}_{1}$ and $\mathbf{P}_{0}\mathbf{R}
$ are collinear (linear dependent), which means mathematically, that 
\begin{equation}
\mathbf{P}_{0}\mathbf{P}_{1}\parallel \mathbf{P}_{0}\mathbf{R:\qquad }%
\left\vert 
\begin{array}{ll}
\left( \mathbf{P}_{0}\mathbf{P}_{1}.\mathbf{P}_{0}\mathbf{P}_{1}\right) & 
\left( \mathbf{P}_{0}\mathbf{P}_{1}.\mathbf{P}_{0}\mathbf{R}\right) \\ 
\left( \mathbf{P}_{0}\mathbf{R}.\mathbf{P}_{0}\mathbf{P}_{1}\right) & \left( 
\mathbf{P}_{0}\mathbf{R}.\mathbf{P}_{0}\mathbf{R}\right)%
\end{array}%
\right\vert =0  \label{a1.14}
\end{equation}

Thus, if the point $P_{1}^{\prime }\in \mathcal{T}_{P_{0}P_{1}}$, $P^{\prime
}\neq P_{1}$, $P^{\prime }\neq P_{0}$, then the straight line $\mathcal{T}%
_{P_{0}P_{1}^{\prime }}=\mathcal{T}_{P_{0}P_{1}}$ in the Euclidean space,
but $\mathcal{T}_{P_{0}P_{1}^{\prime }}\neq \mathcal{T}_{P_{0}P_{1}}$ in the
physical geometry, in general.

Then%
\begin{equation}
\mathcal{CY}\left( P_{0},P_{1},Q\right) =\mathcal{CY}\left(
P_{0},P_{1}^{\prime },Q\right) \wedge P_{1}^{\prime }\in \mathcal{T}%
_{P_{0}P_{1}}  \label{a1.15}
\end{equation}%
in the proper Euclidean geometry. But in general,%
\begin{equation}
\mathcal{CY}\left( P_{0},P_{1},Q\right) \neq \mathcal{CY}\left(
P_{0},P_{1}^{\prime },Q\right) \wedge P_{1}^{\prime }\in \mathcal{T}%
_{P_{0}P_{1}}  \label{a1.16}
\end{equation}%
In other words, the circular cylinder of the proper Euclidean space is split
into many different cylinders after deformation of the Euclidean space.

A more accurate statement is as follows. Cylinders $\mathcal{CY}\left(
P_{0},P_{1},Q\right) $ and $\mathcal{CY}\left( P_{0},P_{1}^{\prime
},Q\right) $ are different, in general. But in the proper Euclidean geometry
they may coincide, even if $P_{1}^{\prime }\neq P_{1}$ but $P_{1}^{\prime
}\in \mathcal{T}_{P_{0}P_{1}}$. It means, that in the proper Euclidean
geometry different geometrical objects may coincide, because of very high
symmetry of the Euclidean geometry. In other words, a deformation of the
Euclidean geometry may violate its symmetry, and coincidence of different
geometrical objects ceases.

\section{Multivariance of two vectors equivalence}

In application to the property of equivalence the multivariance property
looks as follows. In general, there are many vectors $\mathbf{Q}_{0}\mathbf{Q%
}_{1}$,$\mathbf{Q}_{0}\mathbf{Q}_{1}^{\prime },...$ which are equivalent to
vector $\mathbf{P}_{0}\mathbf{P}_{1}$ and are not equivalent between
themselves. This situation is common for all physical geometries. However,
there are possible such geometries, where the set of vectors $\mathbf{Q}_{0}%
\mathbf{Q}_{1}$,$\mathbf{Q}_{0}\mathbf{Q}_{1}^{\prime },...$ degenerates
into one vector. In the proper Euclidean geometry we have such a
degeneration for all vectors $\mathbf{P}_{0}\mathbf{P}_{1}$ and any point $%
Q_{0}$, which is an origin of the vector $\mathbf{Q}_{0}\mathbf{Q}_{1}$,
equivalent to the vector $\mathbf{P}_{0}\mathbf{P}_{1}$.

It means that the multivariance of the equivalence property is a general
property of a physical geometry, whereas single-variance of the equivalence
property in the proper Euclidean geometry is a specific property of the
Euclidean geometry, which is conditioned by the form of the Euclidean world
function.

The multivariance is a new property of physical geometry. Multivariance
properties have not been yet investigated properly. Multivariance of the
equivalence property generates intransitivity of the vector equivalence. In
other words, if $\mathbf{P}_{0}\mathbf{P}_{1}$eqv$\mathbf{Q}_{0}\mathbf{Q}%
_{1}$ and $\mathbf{Q}_{0}\mathbf{Q}_{1}$eqv$\mathbf{S}_{0}\mathbf{S}_{1}$,
then, in general, $\mathbf{P}_{0}\mathbf{P}_{1}$ is not equivalent to $%
\mathbf{S}_{0}\mathbf{S}_{1}$. The intransitive equivalence is difficult for
investigation. It is reasonable to consider physical geometries, which do
not contain multivariance, or contain multivariance in the minimal degree.

The form of the world function is a unique characteristic of a physical
geometry. One can change a physical geometry only changing its world
function. To obtain the Riemannian geometry we are to impose on the world
function the following constraint 
\begin{equation}
F_{3}\left( P_{0},R,P_{1}\right) \equiv \sqrt{2\sigma \left( P_{0},R\right) }%
+\sqrt{2\sigma \left( R,P_{1}\right) }-\sqrt{2\sigma \left(
P_{0},P_{1}\right) }\geq 0,\qquad \forall P_{0},P_{1},R\in \Omega
\label{a1.17}
\end{equation}%
It is the triangle axiom. Its meaning is as follows. The segment $\mathcal{T}%
_{\left[ P_{0}P_{1}\right] }$ of the straight between the points $%
P_{0},P_{1} $ is described by the relation 
\begin{equation}
\mathcal{T}_{\left[ P_{0}P_{1}\right] }=\left\{ R|F_{3}\left(
P_{0},R,P_{1}\right) =0\right\}  \label{a1.18}
\end{equation}%
In general, the equation 
\begin{equation}
F_{3}\left( P_{0},R,P_{1}\right) \equiv \sqrt{2\sigma \left( P_{0},R\right) }%
+\sqrt{2\sigma \left( R,P_{1}\right) }-\sqrt{2\sigma \left(
P_{0},P_{1}\right) }=0  \label{a1.19}
\end{equation}%
describes some surface $S$ around some volume $V$ in the space $\Omega $.
The external points $R$ with respect to the volume $V$ satisfy the relation $%
F_{3}\left( P_{0},R,P_{1}\right) >0$. The internal points $R$ inside the
volume $V$ satisfy the relation $F_{3}\left( P_{0},R,P_{1}\right) <0$. If
the triangle axiom (\ref{a1.17}) is fulfilled, the volume $V$ is empty, and
the segment $\mathcal{T}_{\left[ P_{0}P_{1}\right] }$ of straight is
one-dimensional, because the surface $S$ does not contain any points inside.

On the other hand, the segment (\ref{a1.18}) of the straight can be
described by the relation 
\begin{equation}
\mathcal{T}_{\left[ P_{0}P_{1}\right] }=\left\{ R|\mathbf{P}_{0}\mathbf{P}%
_{1}\upuparrows \mathbf{P}_{0}\mathbf{R\wedge }\left\vert \mathbf{P}_{0}%
\mathbf{R}\right\vert \leq \left\vert \mathbf{P}_{0}\mathbf{P}%
_{1}\right\vert \right\}  \label{a1.20}
\end{equation}%
where $\mathbf{P}_{0}\mathbf{P}_{1}\upuparrows \mathbf{P}_{0}\mathbf{R}$ is
the condition of parallelism of vectors $\mathbf{P}_{0}\mathbf{P}_{1}$ and $%
\mathbf{P}_{0}\mathbf{R}$, which is described by the relation 
\begin{equation}
\mathbf{P}_{0}\mathbf{P}_{1}\upuparrows \mathbf{P}_{0}\mathbf{R:\qquad }%
\left( \mathbf{P}_{0}\mathbf{P}_{1}.\mathbf{P}_{0}\mathbf{R}\right)
=\left\vert \mathbf{P}_{0}\mathbf{R}\right\vert \cdot \left\vert \mathbf{P}%
_{0}\mathbf{P}_{1}\right\vert  \label{a1.21}
\end{equation}

It easy to verify that two definitions (\ref{a1.18}), (\ref{a1.19}) and (\ref%
{a1.20}), (\ref{a1.21}) are equivalent because of the identity%
\begin{equation}
\left( \mathbf{P}_{0}\mathbf{P}_{1}.\mathbf{P}_{0}\mathbf{R}\right)
^{2}-\left\vert \mathbf{P}_{0}\mathbf{R}\right\vert ^{2}\left\vert \mathbf{P}%
_{0}\mathbf{P}_{1}\right\vert ^{2}\equiv \frac{1}{4}F_{0}\left(
P_{0},R,P_{1}\right) F_{1}\left( P_{0},R,P_{1}\right) F_{2}\left(
P_{0},R,P_{1}\right) F_{3}\left( P_{0},R,P_{1}\right)  \label{a1.22}
\end{equation}%
where%
\begin{eqnarray*}
F_{0}\left( P_{0},R,P_{1}\right) &=&\sqrt{2\sigma \left( P_{0},R\right) }+%
\sqrt{2\sigma \left( R,P_{1}\right) }+\sqrt{2\sigma \left(
P_{0},P_{1}\right) } \\
F_{1}\left( P_{0},R,P_{1}\right) &=&\sqrt{2\sigma \left( P_{0},R\right) }-%
\sqrt{2\sigma \left( R,P_{1}\right) }+\sqrt{2\sigma \left(
P_{0},P_{1}\right) } \\
F_{2}\left( P_{0},R,P_{1}\right) &=&-\sqrt{2\sigma \left( P_{0},R\right) }+%
\sqrt{2\sigma \left( R,P_{1}\right) }+\sqrt{2\sigma \left(
P_{0},P_{1}\right) } \\
F_{3}\left( P_{0},R,P_{1}\right) &=&\sqrt{2\sigma \left( P_{0},R\right) }+%
\sqrt{2\sigma \left( R,P_{1}\right) }-\sqrt{2\sigma \left(
P_{0},P_{1}\right) }
\end{eqnarray*}

Thus, if the world function satisfies the triangle axiom (\ref{a1.17}), any
segment (\ref{a1.20}) is single-variant (one-dimensional), and equivalence (%
\ref{a1.3}) of two vectors $\mathbf{P}_{0}\mathbf{P}_{1}$ and $\mathbf{P}_{0}%
\mathbf{R}$ is single-variant, provided they have a common origin $P_{0}$.
It means that it follows from the relation $\mathbf{P}_{0}\mathbf{P}_{1}$eqv$%
\mathbf{P}_{0}\mathbf{R}$ and (\ref{a1.17})\textbf{, }that $R=P_{1}$. Note,
that the single-variance of the two vectors equivalence takes place only for
the proper Riemannian geometry, when the world function is nonnegative and
all terms in the relation (\ref{a1.17}) are real for any points $%
P_{0},P_{1},R\in \Omega $. For the pseudo-Riemannian geometry, when the
world function $\sigma $ may have any sign, the equivalence of two vectors
is multivariant, in general, even if the vectors have a common origin.

Investigation shows , that the equivalence of two vectors is single-variant
only in the flat proper Riemannian space, i.e. in the proper Euclidean
space. Note, that in \cite{R2001} one investigated multivariance of
parallelism of two directions, the concept of the vector equivalence had not
yet been introduced. However, there is a single-valued connection between
the multivariance of two vectors equivalence and multivariance of the two
direction parallelism. There are some special cases of the geometry and of
vectors, when the equivalence property is single-variant. For instance, in
the proper Riemannian space, when the world function is nonnegative and the
triangle axiom (\ref{a1.17}) takes place, equivalence of two vectors $%
\mathbf{P}_{0}\mathbf{P}_{1}$ and $\mathbf{Q}_{0}\mathbf{R}$ is
single-variant, provided that $\mathbf{P}_{0}\mathbf{P}_{1}\parallel \mathbf{%
P}_{0}\mathbf{Q}_{0}$.

\section{Riemannian and $\protect\sigma $-Riemannian geometries}

The physical geometry, whose world function satisfies the triangle axiom (%
\ref{a1.17}) will be referred to as the $\sigma $-Riemannian geometry. The
Riemannian geometry distinguishes from the $\sigma $-Riemannian geometry by
the method of its construction. The $\sigma $-Riemannian geometry is
constructed by means of the deformation principle, and it is defined
completely by its world function. The $\sigma $-Riemannian geometry may be
discrete, or continuous. This circumstance is of no importance, because the
construction of physical geometry by means of the deformation principle does
not need introduction of a coordinate system.

The $n$-dimensional Riemannian geometry is constructed as an internal
geometry of the $n$-dimensional smooth surface $S_{n}$ in $m$-dimensional
proper Euclidean space $E_{m}$ $(m>n)$. One introduces a curvilinear
coordinate system $K_{m}$ in $E_{m}$ with coordinates $\mathbf{\xi }=\left\{ 
\mathbf{x},\mathbf{y}\right\} =\left\{ \xi ^{i}\right\} ,$ $i=1,2,...m$,
where $\mathbf{x}=\left\{ x^{1},x^{2},...x^{n}\right\} $, $\mathbf{y}%
=\left\{ y^{n+1},y^{n+2},...y^{m}\right\} $. Coordinates are introduced in
such a way, that coordinates $\left\{ \mathbf{x},\mathbf{0}\right\} $
describe points of the $n$-dimensional surface $S_{n}$. Besides, coordinates 
$\mathbf{\xi }$ are supposed to be chosen in such a way, that any vector $%
\mathbf{U}\left( P\right) =\left\{ 0,0,...0,U^{n+1},U^{n+2},...U^{m}\right\} 
$ at a point $P\in S_{n}$ is orthogonal to any tangent to the surface $S_{n}$
vector $\mathbf{V}\left( P\right) \mathbf{=}\left\{
V^{1},V^{2},...V^{n},0,0,...0\right\} $, taken at the same point $P\in S_{n}$%
. Let $g_{ik}$, $i,k=1,2,...m$ be the metric tensor in the proper Euclidean
space $E_{m}$ in the coordinate system $K_{m}$. Then%
\begin{equation}
g_{ik}\left( \mathbf{\xi }\right) =g_{ik}\left( \mathbf{x,y}\right)
=\dsum\limits_{l=1}^{l=m}\frac{\partial X^{l}\left( \mathbf{\xi }\right) }{%
\partial \xi ^{i}}\frac{\partial X^{l}\left( \mathbf{\xi }\right) }{\partial
\xi ^{k}},\qquad i,k=1,2,...m,  \label{a4.1}
\end{equation}%
where $X^{l}\left( \mathbf{\xi }\right) ,$ $l=1,2,...m$ are Cartesian
coordinates of the point $\mathbf{\xi }$ in $E_{m}$. Curvilinear coordinates
are chosen in such a way, that on the surface $S_{n}$ they satisfy the
conditions%
\begin{equation}
\dsum\limits_{l=1}^{l=m}\frac{\partial X^{l}\left( \mathbf{x,0}\right) }{%
\partial \xi ^{i}}\frac{\partial X^{l}\left( \mathbf{x,0}\right) }{\partial
\xi ^{k}}=0,\qquad i=1,2,...n,\qquad k=n+1,n+2,...m  \label{a4.2}
\end{equation}

The line element $ds$ on the surface $S_{n}$ is described by the relation%
\begin{equation}
ds^{2}=\dsum\limits_{l=1}^{l=n}g_{ik}\left( \mathbf{x,0}\right)
dx^{i}dx^{k}\equiv g_{ik}\left( \mathbf{x}\right) dx^{i}dx^{k}  \label{a4.3}
\end{equation}%
where the sign of sum is omitted and the summation over repeated indices is
produced from $1$ to $n$. This rule is used further. One can determine
geodesics on the surface $S_{n}$ by means of the relations%
\begin{equation}
\frac{d^{2}x^{k}}{d\tau ^{2}}+\gamma _{ls}^{k}\left( \mathbf{x}\right) \frac{%
dx^{l}}{d\tau }\frac{dx^{s}}{d\tau }=0  \label{a4.4}
\end{equation}%
where $\gamma _{ls}^{k}$ is the Cristoffel symbol in the coordinate system $%
K_{n}$ on the surface $S_{n}$%
\begin{equation}
\gamma _{ls}^{k}\left( \mathbf{x}\right) =\frac{1}{2}g^{ki}\left( \mathbf{x}%
\right) \left( \frac{\partial g_{is}\left( \mathbf{x}\right) }{\partial x^{l}%
}+\frac{\partial g_{li}\left( \mathbf{x}\right) }{\partial x^{s}}-\frac{%
\partial g_{ls}\left( \mathbf{x}\right) }{\partial x^{i}}\right) ,\qquad
k,l,s=1,2,...n  \label{a4.6}
\end{equation}

\[
g^{ki}\left( \mathbf{x}\right) g_{li}\left( \mathbf{x}\right) =\delta
_{l}^{k},\qquad k,l=1,2,...n 
\]%
Let for simplicity there be only one geodesic segment $\mathcal{L}_{\left[
P_{0}P_{1}\right] }\subset S_{n}$, connecting any two points $P_{0},P_{1}\in
S_{n}$. We can define the world function $\sigma _{\mathrm{R}}$ on $S_{n}$
by means of the relation%
\begin{equation}
\sigma _{\mathrm{R}}\left( P_{0},P_{1}\right) =\frac{1}{2}\left(
\dint\limits_{\mathcal{L}_{\left[ P_{0}P_{1}\right] }}\sqrt{g_{ik}\left( 
\mathbf{x}\right) dx^{i}dx^{k}}\right) ^{2},\qquad P_{0},P_{1}\in S_{n}
\label{a4.5}
\end{equation}%
According to definition (\ref{a4.5}) the world function $\sigma _{\mathrm{R}%
} $ satisfies the triangle axiom (\ref{a1.17}). It means, that the world
function $\sigma _{\mathrm{R}}$ of the Riemannian geometry may coincide with
the world function of the $\sigma $-Riemannian geometry, if the set $\Omega $
is identified with the surface $S_{n}$. Then we may repeat construction of
the Riemannian geometry on the surface $S_{n}$ by means of the deformation
principle. In this case we obtain the above obtained results. We obtain the
line element in the form (\ref{a4.3}) and equation for the straight
(geodesic) in the form (\ref{a4.4}). All obtained single-variant results of
the Riemannian geometry may be obtained by means of the deformation
principle from the world function (\ref{a4.5}).

Two vectors $\mathbf{U}\left( \mathbf{x}_{1},\mathbf{0}\right) $ and $%
\mathbf{V}\left( \mathbf{x}_{2},\mathbf{0}\right) $ at two different points $%
P_{1}=\left\{ \mathbf{x}_{1},\mathbf{0}\right\} \in S_{n}$ and $%
P_{2}=\left\{ \mathbf{x}_{2},\mathbf{0}\right\} \in S_{n}$ are equal in the
Euclidean space $E_{m}$, if their Cartesian components coincide%
\begin{equation}
\dsum\limits_{k=1}^{k=m}\frac{\partial X^{l}}{\partial \xi ^{k}}\left( 
\mathbf{x}_{1},\mathbf{0}\right) U^{k}\left( \mathbf{x}_{1},\mathbf{0}%
\right) =\dsum\limits_{k=1}^{k=m}\frac{\partial X^{l}}{\partial \xi ^{k}}%
\left( \mathbf{x}_{2},\mathbf{0}\right) V^{k}\left( \mathbf{x}_{2},\mathbf{0}%
\right) ,\qquad l=1,2,...m  \label{a4.7}
\end{equation}%
where $U^{k}\left( \mathbf{x}_{1},\mathbf{0}\right) $ and $V^{k}\left( 
\mathbf{x}_{2},\mathbf{0}\right) $ are components of vectors $\mathbf{U}%
\left( \mathbf{x}_{1},\mathbf{0}\right) $ and $\mathbf{V}\left( \mathbf{x}%
_{2},\mathbf{0}\right) $ respectively in the curvilinear coordinate system $%
K_{m}$. If vectors $\mathbf{U}\left( \mathbf{x}_{1},\mathbf{0}\right) $ and $%
\mathbf{V}\left( \mathbf{x}_{2},\mathbf{0}\right) $ are vectors of the
internal geometry in $S_{n}$, they are tangent to the surface $S_{n}$, i.e.%
\begin{equation}
U^{k}\left( \mathbf{x}_{1},\mathbf{0}\right) =0,\qquad V^{k}\left( \mathbf{x}%
_{2},\mathbf{0}\right) =0,\qquad k=n+1,n+2,...m  \label{a4.8}
\end{equation}

Let the vector $\mathbf{U}\left( \mathbf{x}_{1},\mathbf{0}\right) $,
satisfying the first relation (\ref{a4.8}) be fixed. Then the vector $%
\mathbf{V}\left( \mathbf{x}_{2},\mathbf{0}\right) $, satisfying conditions (%
\ref{a4.7}), (\ref{a4.8}) does not exist, in general, if the points $P_{1}$
and $P_{2}$ are different, and 
\begin{equation}
\frac{\partial X^{l}}{\partial \xi ^{k}}\left( \mathbf{x}_{1},\mathbf{0}%
\right) \neq \frac{\partial X^{l}}{\partial \xi ^{k}}\left( \mathbf{x}_{2},%
\mathbf{0}\right) ,\qquad l,k=1,2,...m  \label{a4.9}
\end{equation}%
But the multivariant relations of the $\sigma $-Riemannian geometry cannot
be obtained in the Riemannian geometry, because the Riemannian geometry does
not contain multivariant relations in principle. Thus, in general, we cannot
take from the proper Euclidean geometry the concept of the remote vectors
equality in the internal geometry of the surface $S_{n}$.

If the Riemannian geometry is considered as an abstract logical
construction, one may at all not introduce equality of two remote vectors
and their parallelism. In this case we obtain a geometry, which has only
some concepts of the proper Euclidean geometry, but not all of them. The
obtained geometry appears to be more poor in concepts, than the proper
Euclidean geometry.

However, if the (pseudo-)Riemannian geometry pretends to be used for
description of the space-time geometry, one is forced to introduce these
concepts, because a construction of particle dynamics is impossible without
a conception of equality of remote vectors. The parallel transport is
introduced in the Riemannian geometry as follows.

Let $\mathbf{u}\left( \mathbf{x}\right) =\left\{ u^{k}\left( \mathbf{x}%
\right) \right\} ,$ $k=1,2,...n$ be a vector on the surface $S_{n}$.
Simultaneously the vector $\mathbf{U}\left( \mathbf{x,0}\right) =\left\{ 
\mathbf{u}\left( \mathbf{x}\right) ,\mathbf{0}\right\} $ is a vector in the
proper Euclidean space $E_{m}$. Its coordinates in the coordinate system $%
K_{m}$ have the form $U^{k}=u^{k}\left( \mathbf{x}\right) ,$ $k=1,2,...n,$ $%
U^{k}=u^{k}=0,$ $k=n+1,n+2,...m$. Let $d\mathbf{\xi }$ be an infinitesimal
vector of displacement on the surface $S_{n}$, $d\mathbf{\xi =}\left\{ d%
\mathbf{x},\mathbf{0}\right\} $.

Let us transport the vector $\mathbf{U}\left( \mathbf{x,0}\right) $ in $%
E_{m} $ from the point $\mathbf{\xi =}\left\{ \mathbf{x,0}\right\} $ into
the point $\mathbf{\xi }+d\mathbf{\xi =}\left\{ \mathbf{x}+d\mathbf{x,0}%
\right\} $. The transport is produced in such a way, that $\mathbf{U}\left( 
\mathbf{x,0}\right) =\mathbf{U}\left( \mathbf{x+d\mathbf{x},0}\right) $. It
is always possible in the proper Euclidean space $E_{m}$. As far as $d%
\mathbf{x} $ is infinitesimal quantity, in the coordinate system $K_{m}$ the
vector $\mathbf{U}\left( \mathbf{x+d\mathbf{x},0}\right) $ has the form 
\begin{equation}
U^{k}\left( \mathbf{x+d\mathbf{x},0}\right) =U^{k}\left( \mathbf{x,0}\right)
+\delta U^{k}\left( \mathbf{x,0}\right) ,\qquad k=1,2,...m  \label{a4.10}
\end{equation}%
where $\delta U^{k}\left( \mathbf{x,0}\right) $ is an infinitesimal quantity
of the order $O\left( \left\vert d\mathbf{x}\right\vert \right) $. As far as 
$U^{k}\left( \mathbf{x,0}\right) =0$, $k=n+1,n+2,...m$, we obtain 
\begin{equation}
U^{k}\left( \mathbf{x+d\mathbf{x},0}\right) =\delta U^{k}\left( \mathbf{x,0}%
\right) ,\qquad k=n+1,n+2,...m  \label{a4.11}
\end{equation}%
and $\left\vert \mathbf{U}\left( \mathbf{x+d\mathbf{x},0}\right) \right\vert 
$ coincide with $\left\vert \mathbf{U}\left( \mathbf{x,0}\right) \right\vert 
$ to within $O\left( \left\vert d\mathbf{x}\right\vert ^{2}\right) $.

Let us project the vector $\mathbf{U}\left( \mathbf{x+d\mathbf{x},0}\right) $
onto the surface $S_{n}$. It means that we set $\delta U^{k}\left( \mathbf{%
x,0}\right) =0,\qquad k=n+1,n+2,...m$. Calculation of $\delta U^{k}\left( 
\mathbf{x,0}\right) ,\qquad k=1,2,...n$ gives%
\begin{equation}
\delta U^{k}\left( \mathbf{x,0}\right) =\dsum\limits_{l,s=1}^{l,s=n}\gamma
_{ls}^{k}\left( \mathbf{x}\right) U^{l}\left( \mathbf{x,0}\right)
dx^{s}+O\left( \left\vert d\mathbf{x}\right\vert ^{2}\right) ,\qquad
k=1,2,...n  \label{a4.12}
\end{equation}%
where the Christoffel symbol $\gamma _{ls}^{k}\left( \mathbf{x}\right) $ is
given on the surface $S_{n}$ by the relation (\ref{a4.6}).

The relation (\ref{a4.12}) contains only tangential components $\mathbf{u}$
of the vector $\mathbf{U\in E}_{n}$. It means that the relation (\ref{a4.12}%
) is a relation of the internal geometry on the surface $S_{n}$. This
relation may be described in the form, containing only quantities of the
internal geometry. The vector 
\begin{equation}
u^{k}\left( \mathbf{x}+d\mathbf{x}\right) =u^{k}\left( \mathbf{x}\right)
+\delta u^{k}\left( \mathbf{x}\right) =u^{k}\left( \mathbf{x}\right) +\gamma
_{ls}^{k}\left( \mathbf{x}\right) u^{l}\left( \mathbf{x}\right)
dx^{s},\qquad k=1,2,...n  \label{a4.14}
\end{equation}%
is considered to be in parallel with the vector $u^{k}\left( \mathbf{x}%
\right) $, $k=1,2,...n$.

It is well known relation for the parallel transport of a vector in the
Riemannian geometry. The parallel transport from point $\mathbf{x}$ to the
point $\mathbf{x}^{\prime }\mathbf{\ }$is produced as follows. The points $%
\mathbf{x}$ and $\mathbf{x}^{\prime }$ are connected by some line $\mathcal{L%
}_{xx^{\prime }}$. The line is divided into infinitesimal segments. The
parallel transport is produced step by step along all segments by means of
the formula (\ref{a4.14}). Result of the parallel transport depends on the
path $\mathcal{L}_{xx^{\prime }}$. Essentially the result of the parallel
transport is multivariant, but each version of the parallel transport is
connected with some path $\mathcal{L}_{xx^{\prime }}$ of the transport.

The parallel transport in $\sigma $-Riemannian geometry was investigated in
sec.6 of paper \cite{R2001}. Here we present only the result of this
investigation. The parallel transport (parallelism of two vectors) is
defined by the relation (\ref{a1.21}). In the coordinate form it is written
as follows%
\begin{equation}
\left( \sigma _{i,l^{\prime }}\left( x,x^{\prime }\right) \sigma
_{k,s^{\prime }}\left( x,x^{\prime }\right) -g_{ik}\left( x\right)
g_{l^{\prime }s^{\prime }}\left( x^{\prime }\right) \right) u^{i}\left(
x\right) u^{k}\left( x\right) v^{l^{\prime }}\left( x^{\prime }\right)
v^{k^{\prime }}\left( x^{\prime }\right) =0  \label{a4.15}
\end{equation}%
where $\sigma \left( x,x^{\prime }\right) $ is the world function of the $%
\sigma $-Riemannian geometry between the points with coordinates $x=\left\{
x^{i}\right\} $, $i=1,2,...n$ and $x^{\prime }=\left\{ x^{\prime i}\right\} $%
, $i=1,2,...n$ 
\begin{equation}
\sigma _{il^{\prime }}\left( x,x^{\prime }\right) \equiv \sigma
_{i,l^{\prime }}\left( x,x^{\prime }\right) \equiv \frac{\partial ^{2}\sigma
\left( \left( x,x^{\prime }\right) \right) }{\partial x^{i}\partial
x^{\prime l}}  \label{a4.16}
\end{equation}%
Prime at the index means that this index corresponds to the point $x^{\prime
}$. $u^{i}\left( x\right) $ is a contravariant vector at the point $x$. $%
v^{k^{\prime }}\left( x^{\prime }\right) $ is a contravariant vector at the
point $x^{\prime }$. If vectors $u^{k}\left( x\right) $ and $v^{k}\left(
x^{\prime }\right) $ are in parallel, they satisfy the relation (\ref{a4.15}%
). The relation (\ref{a4.15}) has been obtained for infinitesimal vectors $%
u^{i}\left( x\right) $ and $v^{k^{\prime }}\left( x^{\prime }\right) $. But
the relation (\ref{a4.15}) is invariant with respect to a change of lengths
of vectors $u^{i}\left( x\right) $ and $v^{k^{\prime }}\left( x^{\prime
}\right) $, and it appears to be valid also for finite vectors $u^{i}\left(
x\right) $, $v^{k^{\prime }}\left( x^{\prime }\right) $. At the deduction of
the relation (\ref{a4.15}) one uses the fact that the $\sigma $ is the world
function, defined by the relation (\ref{a4.5}), i.e. it is a world function
of the Riemannian space.

Let the vector $\mathbf{x}^{\prime }-\mathbf{x}$ be infinitesimal. Then 
\begin{equation}
d\xi ^{k}\mathbf{=}x^{\prime k}-x^{k}\mathbf{,\qquad }v^{k}\left( \mathbf{x}%
^{\prime }\right) =v^{k}\left( \mathbf{x}\right) +\delta v^{k}\left( \mathbf{%
x}\right) ,\qquad k=1,2,...n  \label{a4.17}
\end{equation}%
where $d\xi ^{k}$ and $\delta v^{k}$ are infinitesimal quantities, and the
relation (\ref{a4.15}) may be transformed to the form%
\begin{eqnarray}
&&\left( u_{k}v^{k}\right) ^{2}-u_{k}u^{k}v_{l}v^{l}+\left( \gamma
_{j;is}g_{kr}+g_{ij}\gamma _{r;ks}-g_{ik}g_{rj,s}\right)
u^{j}u^{r}v^{i}v^{k}d\xi ^{s}  \nonumber \\
&&+\left( g_{ij}g_{kr}-g_{ik}g_{rj}\right) u^{j}u^{r}\delta
v^{i}v^{k}+\left( g_{ij}g_{kr}-g_{ik}g_{rj}\right) u^{j}u^{r}v^{i}\delta
v^{k}=0  \label{a4.18}
\end{eqnarray}%
where 
\begin{equation}
\gamma _{j;is}=g_{jl}\gamma _{is}^{l},\qquad g_{rj,s}\equiv \frac{\partial
g_{rj}}{\partial x^{s}}  \label{a4.19}
\end{equation}%
It is easy to verify, that the conventional parallel transport (\ref{a4.14}) 
\begin{equation}
v^{k}\left( \mathbf{x}\right) =u^{k}\left( \mathbf{x}\right) ,\qquad \delta
v^{k}\left( \mathbf{x}\right) =\gamma _{ls}^{k}\left( \mathbf{x}\right)
v^{l}\left( \mathbf{x}\right) d\xi ^{s}  \label{a4.20}
\end{equation}%
satisfies the relation (\ref{a4.18}). However, the solution (\ref{a4.20}) is
unique, if direction $d\mathbf{\xi }$ of transport coincides with the
direction of the vector $\mathbf{u}\left( \mathbf{x}\right) $. In the
general case, when the curvature of the surface $S_{n}$ does not vanish, the
set of solutions $v+\delta v^{k}$ of equation (\ref{a4.18}) forms a cone
(the collinearity cone). This cone degenerates into one-dimensional line, if 
$\mathbf{u}\left( \mathbf{x}\right) $ is in parallel with $d\mathbf{\xi }$%
\textbf{, }or if the curvature of the surface $S_{n}$ vanishes. In these
cases the conventional parallel transport (\ref{a4.20}) is the unique
solution of (\ref{a4.18}).

One should expect, that the world function determines geometry uniquely. At
any rate the Euclidean world function determines the Euclidean geometry
uniquely. If the world function satisfies the triangle axiom (\ref{a1.17}),
it determines the $\sigma $-Riemannian geometry uniquely. However, there is
in addition the Riemannian geometry which does not coincide with the $\sigma 
$-Riemannian geometry in some details. The Riemannian geometry pretends to
be a true geometry, describing the real space-time. At any rate, most of
contemporary mathematicians consider the Riemannian geometry as a true
geometry and ignore the $\sigma $-Riemannian geometry. We compare the $%
\sigma $-Riemannian geometry and the Riemannian geometry, firstly, as
logical constructions and, secondly, as possible space-time geometries.

The Riemannian geometry is rather special geometry. It is described only in
terms of coordinates. The Riemannian space is isometrically embeddable in
the Euclidean space of rather large dimension. The Riemannian geometry uses
the triangle axiom (\ref{a1.17}), taken in the form (\ref{a4.5}), as
internal constraint of a geometry. It is internal in the sense, that some
basic concepts of the Riemannian geometry (concept of a curve) cannot
formulated at all without the triangle axiom. The concept of the world
function is a secondary concept in the Riemannian geometry. The concept of
the world function, defined by the relation (\ref{a4.5}), cannot be
formulated without a reference to the concept of a curve (geodesic). It is
impossible to obtain a generalization of the Riemannian geometry in terms of
its basic concepts.

The $\sigma $-Riemannian geometry is a special case of the physical
geometry, when the world function is restricted by the triangle axiom (\ref%
{a1.17}). The world function is a primary concept of the $\sigma $%
-Riemannian geometry. The $\sigma $-Riemannian geometry is not constrained
by such conditions as a use of coordinates and isometric embeddability in
the Euclidean space. The $\sigma $-Riemannian geometry does not use such
non-metrical concept as the concept of a curve. The triangle axiom is an
external constraint in the $\sigma $-Riemannian geometry in the sense, that
it is not used in construction of the geometry. Avoiding the triangle axiom,
we obtain a more general geometry, whose concepts are constructed without a
reference to the triangle axiom.

The $\sigma $-Riemannian geometry is a more general construction than, the
Riemannian one in the sense, that imposing some constraints on the $\sigma $%
-Riemannian geometry, one may to obtain the Riemannian geometry. For
instance, in the physical geometry there are two sorts of straights 
\begin{equation}
\mathcal{T}_{P_{0}P_{1};P_{0}}\equiv \mathcal{T}_{P_{0}P_{1}}=\left\{ R|%
\mathbf{P}_{0}\mathbf{P}_{1}\parallel \mathbf{P}_{0}\mathbf{R}\right\}
\label{a4.21}
\end{equation}%
and%
\begin{equation}
\mathcal{T}_{P_{0}P_{1};Q_{0}}=\left\{ R|\mathbf{P}_{0}\mathbf{P}%
_{1}\parallel \mathbf{Q}_{0}\mathbf{R}\right\}  \label{a4.22}
\end{equation}

The straights of the type (\ref{a4.21}) are single-variant (one-dimensional)
in the $\sigma $-Riemannian geometry, whereas the straight of the type (\ref%
{a4.22}) are multivariant, in general. The straight of the type (\ref{a4.21}%
) is a special (degenerate) case of the straight (\ref{a4.22}), when the
point $Q_{0}$ coincides with the point $P_{0}$. In the proper Euclidean
geometry the straights of both types are single-variant (one-dimensional).

To obtain the Riemannian geometry from the $\sigma $-Riemannian geometry,
one needs to remove the multivariant straights of the type (\ref{a4.22}) and
to use only degenerate one-dimensional straights of the type (\ref{a4.21}).
Thereafter the degenerate straights (\ref{a4.21}) are used for introduction
of the concept of a geodesic. The world function is constructed as a
secondary concept on the basis of the infinitesimal line element and of the
geodesic. The world function is not used at the construction of the
Riemannian geometry. The Riemannian geometry is constructed as an internal
geometry of a surface in the proper Euclidean space. Eliminating
multivariant straights (\ref{a4.22}) from the Riemannian geometry and
leaving only degenerate straights, one cannot be sure that the way of
Riemannian geometry construction, based on a use of only degenerate
straights (geodesics), is consistent. A use of geodesics in the construction
of geometry leads to the fact, that two-dimensional Euclidean plane with a
hole cannot be isometrically embedded to the Euclidean plane without a hole.
Besides, in the Riemannian geometry the absolute parallelism is absent,
although it takes place in the $\sigma $-Riemannian geometry.

Thus, it seems, that the Riemannian geometry is only a part of the full
physical geometry. The remaining part of the full geometry is cut by the
constraint, that the degenerate straights (\ref{a4.21}) form the complete
set of straights. This constraint is an internal constraint, which is used
at the construction of the Riemannian geometry. It cannot be removed without
destruction of the Riemannian geometry.

In application to the space-time the $\sigma $-Riemannian geometry is more
effective, than the Riemannian geometry. The $\sigma $-Riemannian geometry
is a more general geometry. To obtain a more general space-time geometry
from the $\sigma $-Riemannian geometry, it sufficient only to remove the
triangle axiom, which is an external constraint. We obtain immediately a
space-time geometry of a general form. In the Riemannian geometry the
triangle axiom is an internal constraint, which is used at the construction
of the Riemannian geometry. Removing the triangle axiom, we lose the method
of the Riemannian geometry construction.

Using the Riemannian geometry as a space-time geometry, one cannot imagine,
that the space-time geometry may be responsible for quantum effects \cite%
{R91}. One cannot imagine that the space-time geometry may be responsible
for a limited divisibility of physical bodies \cite{R2005}. The unlimited
divisibility, used in the Riemannian space-time geometry, generates an
independence of the particle dynamics on the space-time geometry. In
reality, the true space-time geometry must determine the particle dynamics 
\cite{R2007}.

\end{document}